\documentclass[a4paper,10pt]{article}
\usepackage[utf8x]{inputenc}
\usepackage{tracefnt,amsmath,tabu,array}
\usepackage{amssymb,graphicx,setspace,amsfonts,amsbsy}
\usepackage{pifont,latexsym,ifthen,amsthm,rotating,calc,textcase,booktabs,cancel,slashed}

\newtheorem{theorem}{Theorem}[section]
\newtheorem{lemma}[theorem]{Lemma}
\newtheorem{corollary}[theorem]{Corollary}
\newtheorem{proposition}[theorem]{Proposition}
\newtheorem{definition}[theorem]{Definition}

\newtheorem{remark}[theorem]{Remark}
\newcommand{\filledbox}{\leavevmode
  \hbox to.77778em{%
  \hfil\vbox to.675em{\hrule width.6em height.6em}\hfil}}

\newcommand{\Rm}{{\mathbb R}}

\newcommand{\eps}{\varepsilon}

\begin{document}
\tabulinesep=1.0mm
\title{Decay Estimates of High Dimensional Adjoint Radon Transforms}

\author{Ruipeng Shen\\
Centre for Applied Mathematics\\
Tianjin University\\
Tianjin, China
}

\maketitle

\begin{abstract}
 In this paper we prove an optimal $L^2-L^{2d}$ decay estimate of the adjoint Radon transform of compactly supported data in $d$-dimensional space via a geometric method. A similar problem in dimension $3$ has be considered in the author's previous work. This work deals with all higher dimensional case $d\geq 4$. As an application we give the decay of Strichartz norms of $5$-dimensional non-radiative free waves. The general idea is similar to the lower dimensional case but we introduce a new method to prove the corresponding geometric inequality because the old method becomes too complicated in higher dimensions.
\end{abstract}

\section{Introduction} 

\subsection{Background and motivation}

In this work we consider the adjoint Radon transform
\begin{equation} \label{the operator}
 (\mathcal{R}^* G) (x)= \int_{\mathbb{S}^{d-1}} G(x\cdot \omega, \omega) d\omega, \quad x\in \Rm^d. 
\end{equation}
As its name indicates, this operator is exactly the adjoint of the Radon transform defined by (Here $dS$ is the usual measure of the $(d-1)$-dimensional hyperplane $\omega \cdot x = s$.)
\[
 (\mathcal{R} f) (s,\omega) = \int_{\omega\cdot x = s} f(x) dS(x), \qquad (s,\omega) \in \Rm \times \mathbb{S}^{d-1}.
\]
The Radon transform can be applied in partial differential equations and many other aspect of sciences, such as X-ray technology, computer image processing and radio astronomy. Please refer to Helgason \cite{radon1, radonbook} and Ludwig \cite{radon2} for more details about the Radon transforms. 

\paragraph{Application on free waves} We are interested in the application of adjoint Radon transforms on the radiation theory of free waves. The theory of radiation fields was introduced in mathematical physics more than 50 years ago. Please see Friedlander \cite{radiation1, radiation2}, for example. In the past few years this theory has played an important role in the discussion of asymptotic behaviours of non-linear wave equations. We first give a statement of radiation fields. 

\begin{proposition}[Radiation fields, see Duyckaerts-Kenig-Merle \cite{dkm3}] \label{radiation}
Assume that $d\geq 3$ and let $u$ be a solution to the free wave equation $\partial_t^2 u - \Delta u = 0$ with initial data $(u_0,u_1) \in \dot{H}^1 \times L^2(\Rm^d)$. Then ($u_r$ is the derivative in the radial direction)
\[
 \lim_{t\rightarrow \pm \infty} \int_{\Rm^d} \left(|\nabla u(x,t)|^2 - |u_r(x,t)|^2 + \frac{|u(x,t)|^2}{|x|^2}\right) dx = 0
\]
 and there exist two functions $G_\pm \in L^2(\Rm \times \mathbb{S}^{d-1})$ so that
\begin{align*}
 \lim_{t\rightarrow \pm\infty} \int_0^\infty \int_{\mathbb{S}^{d-1}} \left|r^{\frac{d-1}{2}} \partial_t u(r\theta, t) - G_\pm (r\mp t, \theta)\right|^2 d\theta dr &= 0;\\
 \lim_{t\rightarrow \pm\infty} \int_0^\infty \int_{\mathbb{S}^{d-1}} \left|r^{\frac{d-1}{2}} \partial_r u(r\theta, t) \pm G_\pm (r\mp t, \theta)\right|^2 d\theta dr & = 0.
\end{align*}
In addition, the maps $(u_0,u_1) \rightarrow \sqrt{2} G_\pm$ are bijective isometries from $\dot{H}^1 \times L^2(\Rm^d)$ to $L^2 (\Rm \times \mathbb{S}^{d-1})$. 
\end{proposition}

\begin{remark} \label{Hs equivalence}
 The maps $(u_0,u_1) \rightarrow \sqrt{2} G_\pm$ are also bijective isometries from $\dot{H}^\beta \times \dot{H}^{\beta-1}(\Rm^d)$ to $\dot{H}^{\beta-1}(\Rm \times \mathbb{S}^{d-1})$ for $\beta \in \Rm$. Please see Li-Shen-Wang \cite{radiationHs} for its proof. In fact, this result is a direct consequence of the Fourier transform formula of radiation fields given in C\^{o}te and Laurent's recent work \cite{newradiation}.
\end{remark}

\paragraph{Explicit formula} In this work we call the functions $G_\pm$ radiation profiles. The map from radiation profile $G_-(s,\omega)$ to the corresponding free wave $u$ can be given explicitly 
\begin{equation} \label{from radiation to free wave}
 u(x,t) = \frac{1}{(2\pi)^\mu} \int_{\mathbb{S}^{d-1}}  G_-^{(\mu-1)} \left(x \cdot \omega +t , \omega\right) d\omega.
\end{equation}
Here $\mu = (d-1)/2$ is a constant depending on the dimension $d$ and $G_-^{(\mu-1)}$ is the $(\mu-1)$th partial derivative of $G_-$ with respect to the first derivative $s$. If $d$ is even, then $\mu$ is an half integer. In this case we have to utilize the half derivative operator. The details of this formula can be found, for example, in Li-Shen-Wei \cite{shenradiation}. We may rewrite this formula in the form of the adjoint Radon transform
\[
 u(\cdot,t) = \mathcal{R}^* G_{\mu-1, t} 
\]
Here $G_{\mu-1,t}$ is a time-translated version of the derivative of $G$ 
\[
 G_{\mu-1,t} (s,\omega)= G^{(\mu-1)} (s+t, \omega). 
\]

\paragraph{Non-radiative solutions} We are particularly interested in the case when $d$ is odd and $G_-$ is compactly supported in $[-R,R]\times \mathbb{S}^{d-1}$. This is equivalent to saying that the corresponding free wave $u$ is $R$-weakly non-radiative, i.e.
\[
 \lim_{t\rightarrow \pm \infty} \int_{|x|>|t|+R} |\nabla_{t,x} u|^2 dx = 0.
\]
The non-radiative solutions play an essential role in the channel of energy method, which has many important applications in the study of non-linear wave equations, such as the soliton resolution of solutions to the focusing, energy-critical wave equation (see Duyckaerts-Kenig-Merle \cite{se, oddhigh}) and conditional scattering of solutions to energy super or sub-critical wave equations (see, for instance, Duyckaerts-Kenig-Merle \cite{dkm2} and Shen \cite{shen2}). 

\paragraph{Goal} In this work we will give a decay estimate concerning the adjoint Radon transform of a function $G$ whose support is contained in $[-b,b]\times \mathbb{S}^{d-1}$:
\[
 \|\mathbf{T} G\|_{L^{2d}(\{x: |x|>R\})} \lesssim (b/R)^{-\frac{d-1}{2d}} \|G\|_{L^2}. 
\]
This gives a power-type decay of the Strichartz norm of non-radiative solutions (we give the $5$-dimensional case for an example)
\[
 \|u\|_{L^{7/3} L^{14/3}(\Rm \times \{x: |x|>r\})} \lesssim r^{-2/35}\|(u_0,u_1)\|_{\dot{H}^1\times L^2(\Rm^5)}. 
\]
Following the same argument as in Li et al.\cite{nonradialCE}, we may use this decay estimate of non-radiative solutions to show that the non-radiative solutions to a wide range of energy critical wave equations share the same asymptotic behaviours as the non-radiative free waves. This kind of similarity helps to study asymptotic behaviours of non-linear solutions (see, for example, Duyckaerts-Kenig-Merle \cite{oddtool, oddhigh}) and is an important aspect of the channel of energy method. Most of previous results of this kind apply to radial solutions only.  This paper (as well as \cite{inequality26} in dimension $3$) is an attempt to generalize this theory to non-radial solutions. 

\subsection{Main idea} 

In this work we follow the same idea as given in our previous work \cite{inequality26}, which proves the decay estimates of the three-dimensional adjoint Radon transform. The argument consists of two major steps
\begin{itemize}
 \item[(I)] We first reduce the problem to an integral inequality concerning reciprocal geometric objects, by applying a Cauchy-Schwarz inequality. In the case of three-dimensional adjoint Radon transform we consider an integral over all reciprocal triangles of an arbitrary fixed triangle with vertices in an annulus. The conception of reciprocal objects will be discussed in details later. 
 \item[(II)] We then prove the integral inequality via a few geometric observations. In order to deal with reciprocal triangles, we first classifies them by their sizes and shapes, then discusses their relative locations with respect to the fixed triangle. 
\end{itemize}
It turns out that the original argument in the first step applies to all dimensions with minor modifications. In the second step, however, the original argument could become more complicated and even unrealistic as the dimension grows. In this work we give a new argument to deal with the high dimensional case, which depends on a change of variables method and a relatively simple geometric observation. We first recall the conception of reciprocal triangles and give the definition of high-dimensional reciprocal simplexes.

\paragraph{The reciprocal conception} We recall the conception of reciprocal triangles (or triples) introduced in Li-Shen-Wang \cite{inequality26}. Given six points $A_1, A_2, \cdots, A_6$ in the plane $\Rm^2$, there are many different ways to split them into two group of three points, or equivalently, two triangles. If the product of the areas of these two triangles takes a maximum among all different grouping method, we call these two triples (or triangles) reciprocal to each other. Namely, $A_1 A_2 A_3$ and $A_4 A_5 A_6$ are reciprocal if and only if 
\[
 |A_1 A_2 A_3| \cdot |A_4 A_5 A_6| = \max_{\{j_1,j_2,\cdots, j_6\}=\{1,2,\cdots, 6\}} |A_{j_1} A_{j_2} A_{j_3}| \cdot |A_{j_4} A_{j_5} A_{j_6}|.
\]
These can be generalized to higher dimensional space

\begin{definition}
 We call two simplxes $A_1 A_2 \cdots A_{d+1}$ and $A_{d+2} A_{d+3} \cdots A_{2d+2}$ in $\Rm^d$ reciprocal to each other, if and only if the product of their volumes satisfies 
 \[
  |A_1 A_2 \cdots A_{d+1}|\cdot |A_{d+2} A_{d+3} \cdots A_{2d+2}| = \max |A_{j_1} A_{j_2} \cdots A_{j_{d+1}}|\cdot |A_{j_{d+2}} A_{j_{d+3}} \cdots A_{j_{2d+2}}|.
 \]
 Here the maximum is taken for all permutations $\{j_1,j_2,\cdots,j_{2d+2}\}$ of the first $2d+2$ positive integers $1,2,\cdots, 2d+2$. 
\end{definition}

\subsection{Main results}

In this subsection we give the main results of this work and a few corollaries and remarks. There are three main results. The first one is the geometric inequality; the second one is the decay estimate of the adjoint Radon transforms; the last one is the Strichartz decay of non-radiative free waves, as an application of the decay of adjoint Radon transforms. 

\begin{theorem} \label{main 1}
 Assume $d\geq 3$. Let $\Omega_\ast$ be a sphere shell of outer radius $r$ and thickness $w$ and $A_0 A_1 \cdots A_d$ be an arbitrary fixed simplex in $\Rm^d$. Then we have 
 \[
  \int_{\Sigma(A_0 A_1 \cdots A_d)\cap \Omega_\ast^{d+1}} \frac{dX_0 dX_1 \cdots dX_{d}}{|X_0 X_1 \cdots X_d|} \lesssim_d w^{d+1} r^{d^2-d-1}.
 \]
 Here $\Sigma(A_0 A_1 \cdots A_d)$ is the region in $(\Rm^d)^{d+1}$ consisting of all reciprocal simplexes of $A_0 A_1 \cdots A_d$. Thus the integral region consists of all reciprocal simplexes of $A_0 A_1 \cdots A_d$ with vertices in the sphere shell $\Omega_\ast$. 
\end{theorem}

\begin{remark}
 The upper bound $w^{d+1} r^{d^2-d-1}$ given in Theorem \ref{main 1} is optimal (up to an absolute constant). Please see Subsection \ref{sec: optimal upper bound} for more details. 
\end{remark}

\begin{theorem} \label{main 2}
 Let $d\geq 4$. The adjoint Radon operator $\mathcal{R}^*$ satisfies the following decay estimates
 \begin{itemize}
  \item[(a)] Assume that $G \in L^2(\Rm\times \mathbb{S}^{d-1})$ is supported in $[a,b] \times \mathbb{S}^{d-1}$. Here the constant $b>a>0$ satisfies $b/a<\frac{d^2-2d-1}{2(d-1)}$. Then 
   \[
    \int_{|x|>R} \left|(\mathcal{R}^* G)(x)\right|^{2d} dx \lesssim_d \frac{(b-a)^d}{a (\max\{R, a\})^{d-1}} \|G\|_{L^2}^{2d}.
   \]
   In particular, we have the global estimate
   \[
    \int_{\Rm^d} \left|(\mathcal{R}^* G)(x)\right|^{2d} dx \lesssim_d (b/a-1)^d \|G\|_{L^2}^{2d}.
   \]
  \item[(b)] If $G \in L^2(\Rm\times \mathbb{S}^{d-1})$ is supported in $[-b,b] \times \mathbb{S}^{d-1}$ with $b>0$, then 
  \[
   \int_{|x|>R} \left|(\mathcal{R}^* G)(x)\right|^{2d} dx \lesssim_d (b/R)^{d-1} \|G\|_{L^2}^{2d}, \qquad R\geq b. 
  \]
 \end{itemize}
\end{theorem}

 \begin{remark}
 The decay rate $\|\mathcal{R}^* G\|_{L^{2d} (\{x: |x|>R\})} \lesssim R^{-\frac{d-1}{2d}} \|G\|_{L^2}$ given above is optimal for large $R$'s. We define ($\omega = (\omega_1, \omega_2,\cdots, \omega_d) \in \mathbb{S}^{d-1}$)
\[
 G(s,\omega) = \left\{\begin{array}{ll} R^{1/2}, & \hbox{if}\; s\in [-1,1], \; 0<\omega_d<\frac{1}{4R};\\
 0, & \hbox{otherwise;}  \end{array}\right.
\] 
A basic calculation shows that $\|G\|_{L^2} \simeq 1$. In addition, we consider the value of $\mathcal{R}^* G (x)$ when $x$ is in a circular cylinder $C_R \doteq \{x: |x_d|\leq 2R, x_1^2 + x_2^2 + \cdots + x_{d-1}^2 < 1/4\}$. If $x\in C_R$ and $\omega \in \mathbb{S}^{d-1}$ satisfies $|\omega_d| < 1/4R$, then we have $|x\cdot \omega|<1$. Thus 
\[
 \mathcal{R}^* G (x) = \int_{\omega \in \mathbb{S}^{d-1}, 0<\omega_d < 1/4R} G(x\cdot \omega,\omega) d\omega \simeq R^{-1/2}, \qquad  x\in C_R. 
\]
Therefore we have $\|\mathbf{T} G\|_{L^{2d}(\{x: |x|>R\})} \gtrsim R^{-1/2}\cdot R^{1/2d} \simeq R^{-\frac{d-1}{2d}} \|G\|_{L^2}$. 
\end{remark}

\begin{corollary} \label{global layer like Lp}
 Given $\gamma > 1$, the adjoint Radon transform satisfies the following estimate
 \[
  \left(\sum_{k=-\infty}^\infty \|\mathcal{R}^* G\|_{L^{2d}(\{x: \gamma^k < |x|\leq \gamma^{k+1}\})}^2\right)^{1/2} \lesssim_{d,\gamma} \|G\|_{L^2(\Rm \times \mathbb{S}^{d-1})}
 \]
\end{corollary} 

\begin{remark}
 It has been known that the Radon transform $\mathcal{R}$ is a bounded operator from $L^{2d/(2d-1)}(\Rm^d)$ to $L^2(\Rm\times \mathbb{S}^{d-1})$. This is a special case of the $L^p-L^q$ type estimates of the Radon transform given in Oberlin-Stein \cite{steinRadon}. It immediately follows that the joint Radon transform must be a bounded operator from $L^2 (\Rm \times \mathbb{S}^{d-1})$ to $L^{2d}(\Rm^d)$. Corollary \ref{global layer like Lp} is actually a slightly stronger version of this $L^2-L^{2d}$ estimate since we always have 
 \begin{align*}
 \|\mathcal{R}^* G\|_{L^{2d}(\Rm^d)} & = \left(\sum_{k=-\infty}^\infty \|\mathcal{R}^* G\|_{L^{2d}(\{x: \gamma^k < |x|\leq \gamma^{k+1}\})}^{2d}\right)^{1/2d}\\
 &\qquad \leq \left(\sum_{k=-\infty}^\infty \|\mathcal{R}^* G\|_{L^{2d}(\{x: \gamma^k < |x|\leq \gamma^{k+1}\})}^2\right)^{1/2}.
 \end{align*}
 By the basic properties of adjoint operators, we also have the following estimate of Radon transforms: if $f \in L^{2d/(2d-1)} (\Rm^d)$ is always zero for $|x|<R$, then 
 \[
  \|\mathcal{R} f\|_{L^2([-r,r]\times \mathbb{S}^{d-1})} \lesssim_d (r/R)^{-\frac{d-1}{2d}} \|f\|_{L^{2d/(2d-1)} (\Rm^d)}.
 \]
\end{remark}

\begin{corollary} \label{translation narrow}
 Let $d\geq 4$. If $G \in L^2(\Rm\times \mathbb{S}^{d-1})$ is supported in $[a,b] \times \mathbb{S}^{d-1}$, then for any $R>0$, we have 
 \[
  \int_{|x|>R} \left|(\mathcal{R}^* G)(x)\right|^{2d} dx \lesssim_d \left(\frac{b-a}{R}\right)^{d-1} \|G\|_{L^2}^{2d}
 \]
\end{corollary}

\begin{proposition}
 If $u$ is an $R$-weakly non-radiative free wave in the $5$-dimensional space with a finite energy $E$, then 
 \[
  \|u\|_{L^{7/3} L^{14/3}(\Rm \times \{x\in \Rm^5: |x|>r\})} \lesssim (R/r)^{2/35} E^{1/2}.
 \]
\end{proposition}

\section{Geometric Inequality}

In this section we prove our first main result, i.e. Theorem \ref{main 1}. We start by a few geometric observations, then give a change of variables formula, finally we combine the observations with change of variables to prove Theorem \ref{main 1}.

\subsection{Geometric Observations}

\begin{lemma} \label{location of reciprocal}
 Let $A_0 A_1 \cdots A_d$ and $B_0 B_1 \cdots B_d$ be reciprocal simplexes in $\Rm^d$. Then we have 
 \[
  \min_{k\in\{0,1,\cdots,d\}} \hbox{dist} (B_k, A_0 A_1 \cdots A_{d-1}) \lesssim \hbox{dist} (A_d, A_0 A_1 \cdots A_{d-1})
 \]
 Here $\hbox{dist} (X, A_0 A_1 \cdots A_{d-1})$ represents the distance from the point $X\in \Rm^d$ to the hyperplane $A_0 A_1 \cdots A_{d-1}$.
\end{lemma}
\begin{proof}
 It is clear that there exists $k \in \{0,1,\cdots,d\}$, so that 
\[
  |A_d B_0 B_1 \cdots B_{k-1} B_{k+1} \cdots B_d| \geq \frac{1}{d+1} |B_0 B_1 \cdots B_d|. 
\]
The reciprocal assumption then gives 
\[
  |B_k A_0 A_1 \cdots A_{d-1}| \lesssim |A_0 A_1 \cdots A_d|,
\]
This finishes the proof.
\end{proof}

\begin{lemma} \label{intersection sphere shell}
 Let $\Omega=\{x\in \Rm^d: r-w \leq |x| \leq r\}$ be a sphere shell of outer radius $r$ and thickness $w$. Then the intersection $\Omega_\ast$ of $\Omega$ with any $d-1$-dimensional hyperplane is a $(d-1)$-dimensional sphere shell (or a sphere). In addition, the outer radius $r_\ast$ and thickness $w_\ast$ of $\Omega_\ast$ satisfy the inequality $r_\ast w_\ast \leq 2rw$. Please note that a sphere can be viewed as a sphere shell whose thickness and outer radius are the same. 
\end{lemma}
\begin{proof}
 This lemma follows a direct computation. By radial symmetry it suffices to consider the hyperplane $x_d = c$. The intersection can be given by
\[
 \{x: x_d = c, \; (r-w)^2 - c^2 \leq x_1^2 + x_2^2 + \cdots + x_{d-1}^2 \leq r^2 - c^2\}
\]
This is either a sphere shell, if $|c|<r-w$; a sphere, if $r-w \leq |c| < r$; a single point, if $|c|=r$; or none, if $|c|>r$. In the sphere shell case, we have 
\begin{align*}
 r_\ast w_\ast & = \sqrt{r^2-c^2} \left(\sqrt{r^2-c^2} - \sqrt{(r-w)^2 - c^2}\right) = \frac{\sqrt{r^2-c^2}}{\sqrt{r^2-c^2} + \sqrt{(r-w)^2 - c^2}} w(2r-w)\\
 & \leq 2rw.
\end{align*}
In the sphere case, we have 
\[
 r_\ast w_\ast = r_\ast^2 = r^2-c^2 \leq r^2 - (r-w)^2 \leq 2wr.
\]
\end{proof}

\subsection{Change of variables} \label{sec: change of variables}

Let $A_0, A_1, \cdots, A_d$ be $d+1$ points in the space $\Rm^d$ so that they are not contained in the same $(d-1)$-dimensional hyperplane. Let $\mathbf{x}^{(k)} = (x_1^{(k)}, x_2^{(k)}, \cdots, x_d^{(k)})$ be the coordinates of the point $A_k$. We now introduce a new set of variables. First of all, we let 
\[
 \mathbf{y}^{(0)} = (y_1^{(0)}, y_2^{0}, \cdots, y_d^{(0)}) = (x_1^{(0)}, x_2^{(0)}, \cdots, x_d^{(0)}),
\]
and 
\[
 \left(\prod_{j=1}^{d-1} \sin \omega_j,\; \cos \omega_{d-1} \prod_{j=1}^{d-2} \sin \omega_j,\; \cos \omega_{d-2} \prod_{j=1}^{d-3} \sin \omega_j,\; \cdots,\; \cos \omega_1 \right) \in \mathbb{S}^{d-1} 
\]
be the unit normal vector of the hyperplane $A_0 A_1 \cdots A_{d-1}$ pointing to the direction of the point $A_d$. Here 
\[
  \omega = (\omega_1, \omega_2, \cdots, \omega_{d-1}) \in [0,\pi] \times [0,\pi] \times \cdots \times [0,\pi] \times [0,2\pi)
\]
We may define an orthogonal matrix associated to the unit vector given above 
\[
 S = \begin{pmatrix}
   \displaystyle \prod_{j=1}^{d-1} \sin \omega_j & \cos \omega_{d-1} & \cos \omega_{d-2}  \sin \omega_{d-1} & \cdots & \cdots & \displaystyle \cos \omega_1 \prod_{j=2}^{d-1} \sin \omega_j\\
   \displaystyle \cos \omega_{d-1} \prod_{j=1}^{d-2} \sin \omega_j & -\sin \omega_{d-1} &  \cos \omega_{d-2} \cos \omega_{d-1} & \cdots & \cdots & \displaystyle \cos \omega_1 \cos \omega_{d-1} \prod_{j=2}^{d-2} \sin \omega_j \\
   \displaystyle \cos \omega_{d-2} \prod_{j=1}^{d-3} \sin \omega_j & 0 & -\sin \omega_{d-2} & \cdots & \cdots & \displaystyle \cos \omega_1 \cos \omega_{d-2} \prod_{j=2}^{d-3} \sin \omega_j \\
   \displaystyle \cos \omega_{d-3} \prod_{j=1}^{d-4} \sin \omega_j & 0 & 0 & \cdots & \cdots & \displaystyle \cos \omega_1 \cos \omega_{d-3} \prod_{j=2}^{d-4} \sin \omega_j\\
   \vdots & \vdots & \vdots & \ddots & & \vdots\\
   \vdots & \vdots & \vdots & & \ddots & \vdots \\
   \cos \omega_2 \sin \omega_1 & 0 & 0 & \cdots & \cdots & \cos \omega_1 \cos \omega_2\\
   & & & & & \\
   \cos \omega_1& 0 & 0 & \cdots & \cdots & -\sin \omega_1
 \end{pmatrix}
\]
Here the first column $\Theta_0$ is exactly the normal vector we defined above. We also use the notation $\Theta_{k}$ for the $(k+1)$th column vector in the matrix above. For the reader's convenience we also give the coordinates of $\Theta_k = (\Theta_k^1, \Theta_k^2, \cdots \Theta_k^d)^T$
\[
 \Theta_{k}^m = \left\{\begin{array}{ll} \displaystyle \cos \omega_{d-k}\prod_{j=d+1-k}^{d-1} \sin \omega_j, & m=1;\\
 \displaystyle \cos \omega_{d-k} \cos \omega_{d+1-m} \prod_{j=d+1-k}^{d-m} \sin \omega_j, & 2\leq m \leq k;\\
 -\sin \omega_{d-k}, & m = k+1;\\
 0, & m>k+1.
  \end{array} \right.
\]
\begin{figure}[h]
 \centering
 \includegraphics[scale=1.25]{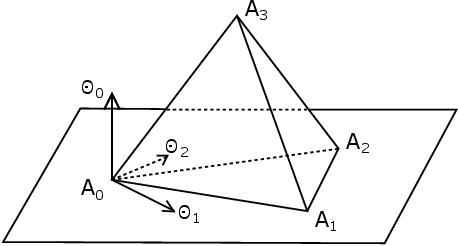}
 \caption{New coordinate system} \label{figure newcoordinates}
\end{figure}
The vectors $\Theta_j$ immediately give a new coordinate systems in $\Rm^d$. More precisely we may choose $A_0$ to be origin of the new system and let $\{\Theta_j\}_{j=0,1,\cdots, d-1}$ be the direction vectors of the new coordinate axes. Figure \ref{figure newcoordinates} gives an illustration of the new coordinate system(here we give the 3-dimensional case as an example). We then let $\mathbf{y}^{(d)} = (y_0^{(k)}, y_1^{(k)}, \cdots, y_{d-1}^{(k)})$ be the coordinates of $A_d$ in the new coordinate system. This is equivalent to saying 
\[
 \overrightarrow{A_0 A_{d}} = \sum_{j=0}^{d-1} y_j^{(d)} \Theta_j,
\]
Similarly we may construct a new coordinate system in the hyperplane $A_0 A_1 \cdots A_{d-1}$. In fact the hyperplane $A_0 A_1 \cdots A_{d-1}$ becomes a coordinate hyperplane in the new coordinate system above. This naturally gives $(d-1)$ coordinates $\mathbf{y} = (y_1, y_2, \cdots, y_{d-1})$ for each point in the hyperplane. Given $k \in \{1,2,\cdots, d-1\}$, we use the notation $\mathbf{y}^{(k)} = (y_1^{(k)}, y_2^{(k)}, \cdots, y_{d-1}^{(k)})$ for the new coordinates of $A_k$ in the hyperplane. We may also write
\[
 \overrightarrow{A_0 A_k} = \sum_{j=1}^{d-1} y_j^{(k)} \Theta_j, \qquad k=1,2,\cdots, d-1.
\]
In later part of this section we will use the change of variable 
\[
 \mathbf{T}: \left(\mathbf{x}^{(0)}, \mathbf{x}^{(1)}, \cdots, \mathbf{x}^{(d)}\right) \rightarrow \left(\omega, \mathbf{y}^{(0)}, \mathbf{y}^{(1)}, \cdots, \mathbf{y}^{(d)}\right).
\]
Next we introduce the change of variable formula 
\begin{lemma}
 Let $\{\mathbf{x}^{(k)}\}$ and $\{\omega,\mathbf{y}^{(k)}\}$ be variables representing points $A_0 A_1 \cdots A_d$ as given above. Then we have the change of variables formula 
\[
 d\mathbf{x}^{(0)} \cdots d\mathbf{x}^{(d)} = |A_0 A_1 \cdots A_{d-1}| d\mathbb{S}^{d-1} (\omega) d\mathbf{y}^{(0)} \cdots d\mathbf{y}^{(d)}.
\]
Here $|A_0 A_1 \cdots A_{d-1}|$ is the surface area of $A_0 A_1\cdots A_{d-1}$ in the $d-1$-dimensional hyperplane; $d\mathbb{S}^{d-1}$ uses the standard surface measure of the $(d-1)$-dimensional sphere $\mathbb{S}^{d-1}$.
\end{lemma}
\begin{proof}
 It suffice to calculate the determinant of Jacobi matrix 
 \[
 J =  \begin{vmatrix}
   M_{0,0} & M_{0,1} & M_{0,2} & \cdots & M_{0,d-1} & M_{0,\omega} & M_{0,d} \\
   M_{1,0} & M_{1,1} & M_{1,2} & \cdots & M_{1,d-1} & M_{1,\omega} & M_{1,d} \\
   \vdots & \vdots & \vdots & \ddots & \vdots & \vdots & \vdots\\
   M_{d,0} & M_{d,1} & M_{d,2} & \cdots & M_{d,d-1} & M_{d,\omega} & M_{d,d}
  \end{vmatrix}
 \]
 Here $M_{k, k'}$ and $M_{k, \omega}$ are all ``local'' Jacobi matrices 
 \begin{align*}
  &M_{k,k'} = \frac{\partial \mathbf{x}^{(k)}}{\partial \mathbf{y}^{(k')}};& &M_{k,\omega} = \frac{\partial \mathbf{x}^{(k)}}{\partial \omega}.&
 \end{align*}
 More precisely we have 
 \begin{align*}
  &M_{k,0} = \frac{\partial(x_1^{(k)}, x_2^{(k)}, \cdots, x_d^{(k)})}{\partial(y_1^{(0)}, y_2^{(0)}, \cdots, y_d^{(0)})};& &M_{k,d} = \frac{\partial(x_1^{(k)}, x_2^{(k)}, \cdots, x_d^{(k)})}{\partial(y_0^{(d)}, y_1^{(d)}, \cdots, y_{d-1}^{(d)})};& \\
  &M_{k,k'} = \frac{\partial(x_1^{(k)}, x_2^{(k)}, \cdots, x_d^{(k)})}{\partial(y_{1}^{(k')}, y_2^{(k')}, \cdots, y_{d-1}^{(k')})},& &k' \in \{1,2,\cdots, d-1\};& \\
  &M_{k,\omega}  = \frac{\partial(x_1^{(k)}, x_2^{(k)}, \cdots, x_d^{(k)})}{\partial(\omega_1, \omega_2, \cdots, \omega_{d-1})}.& &&
\end{align*}
In order to calculate the partial derivatives, we observe
 \begin{align*}
  &\mathbf{x}^{(0)} = \mathbf{y}^{(0)};& &\mathbf{x}^{(d)} = \mathbf{y}^{(0)}+\sum_{j=0}^{d-1} y_j^{(d)}\Theta_j; & &\mathbf{x}^{(k)} = \mathbf{y}^{(0)} + \sum_{j=1}^{d-1} y_j^{(k)}\Theta_j, \qquad k=1,2,\cdots, d-1.&
 \end{align*}
 It immediately follows that 
 \begin{align*}
  &M_{0,k} = M_{0,\omega} = 0, \quad k\geq 1; & &M_{k,d} = 0, \qquad k<d.&
 \end{align*}
 Therefore we have 
 \[
  J = |M_{0,0}|\cdot |M_{d,d}| \cdot \begin{vmatrix}
   M_{1,1} & M_{1,2} & \cdots & M_{1,d-1} & M_{1,\omega} \\
   M_{2,1} & M_{2,2} & \cdots & M_{2,d-1} & M_{2,\omega} \\
       \vdots & \vdots & \ddots & \vdots & \vdots \\
   M_{d-1,1} & M_{d-1,2} & \cdots & M_{d-1,d-1} & M_{d-1,\omega} 
  \end{vmatrix}
 \]
 A basic calculation shows that $M_{0,0}$ is the $d$ by $d$ identity matrix and $M_{d,d} = S$ is an orthogonal matrix. In addition, we have 
 \begin{align*}
  &M_{k,k'} = 0, \qquad k\neq k', \; 1\leq k,k'\leq d-1; & &M_{k,k} = \mathbf{\Theta}, \qquad 1\leq k \leq d-1.&
 \end{align*}
 Here $\mathbf{\Theta} = (\Theta_1, \Theta_2, \cdots, \Theta_{d-1})$ is a $d\times(d-1)$ matrix whose $j$-th columns are exactly the unit vectors $\Theta_j$. Therefore we have 
 \[
  J = \begin{vmatrix}
   \mathbf{\Theta} & 0 & \cdots & 0 & M_{1,\omega} \\
   0 & \mathbf{\Theta} & \cdots & 0 & M_{2,\omega} \\
       \vdots & \vdots & \ddots & \vdots & \vdots \\
   0 & 0 & \cdots & \mathbf{\Theta} & M_{d-1,\omega} 
  \end{vmatrix}
 \]
 We may calculate the $l$-th column $M_{k,\omega}^l$ of the matrix $M_{k,\omega}$ 
 \[
  M_{k,\omega}^l = \sum_{j=1}^{d-1} y_j^{(k)} \frac{\partial \Theta_j}{\partial \omega_l} = \sum_{i=0}^{d-1} \left( \sum_{j=1}^{d-1} y_j^{(k)} \frac{\partial \Theta_j}{\partial \omega_l} \cdot \Theta_i\right) \Theta_i.
 \]
 We may apply basic column operations in the square matrix above and obtain 
 \[
  J = \begin{vmatrix}
   \mathbf{\Theta} & 0 & \cdots & 0 & v_{1,1}\Theta_0 & v_{1,2}\Theta_0 & \cdots & v_{1,d-1} \Theta_0  \\
   0 & \mathbf{\Theta} & \cdots & 0 & v_{2,1}\Theta_0 & v_{2,2}\Theta_0 & \cdots & v_{2,d-1} \Theta_0 \\
       \vdots & \vdots & \ddots & \vdots & \vdots & \vdots & \ddots& \vdots \\
   0 & 0 & \cdots & \mathbf{\Theta} & v_{d-1,1}\Theta_0 & v_{d-1,2}\Theta_0 & \cdots & v_{d-1,d-1} \Theta_0  
  \end{vmatrix}
 \]
Here $v_{k,l}$ is defined by 
\[
 v_{k,l} =  \sum_{j=1}^{d-1} y_j^{(k)} \frac{\partial \Theta_j}{\partial \omega_l} \cdot \Theta_0.
\]
Next we may transform the matrix $\{v_{k,l}\}$ into an upper triangle matrix by switching its columns and/or adding a multiple of one column to another one:
\[
 \begin{pmatrix} v_{1,1} & v_{1,2}& \cdots & v_{1,d-1} \\ v_{2,1} & v_{2,2} & \cdots & v_{2,d-1}\\ \vdots & \vdots & \ddots & \vdots\\ v_{d-1,1} & v_{d-1,2} & \cdots & v_{d-1,d-1}\end{pmatrix} \rightarrow \begin{pmatrix} c_{1,1} & c_{1,2} & \cdots & c_{1,d-1} \\ 0 & c_{2,2} & \cdots & c_{2,d-1}\\ \vdots & \vdots & \ddots & \vdots\\ 0 & 0 & \cdots & c_{d-1,d-1}\end{pmatrix}
\] 
Applying the corresponding column transformations in the last $(d-1)$ columns of the big matrix above, we obtain 
\[
  J = \pm \begin{vmatrix}
   \mathbf{\Theta} & 0 & \cdots & 0 & c_{1,1} \Theta_0 & c_{1,2} \Theta_0 & \cdots& c_{1,d-1} \Theta_0  \\
   0 & \mathbf{\Theta} & \cdots & 0 & 0 & c_{2,2} \Theta_0 & \cdots&  c_{2, d-1} \Theta_0 \\
       \vdots & \vdots & \ddots & \vdots & \vdots & \vdots & \ddots & \vdots \\
   0 & 0 & \cdots & \mathbf{\Theta} & 0 & 0 & \cdots & c_{d-1,d-1} \Theta_0
  \end{vmatrix}
 \]
We may further switching the columns and obtain 
\[
 J = \pm \begin{vmatrix}
   \mathbf{\Theta}_1^+ & \ast & \cdots & \ast\\
   0 & \mathbf{\Theta}_2^+ & \cdots & \ast \\
       \vdots & \vdots & \ddots & \vdots \\
   0 & 0 & \cdots & \mathbf{\Theta}_{d-1}^+ 
  \end{vmatrix}
\]
Here  $\mathbf{\Theta}_k^+ = \left(\mathbf{\Theta}\; c_{k,k}\Theta_0\right)$ are all $d\times d$ matrices. We recall the orthogonality of $\Theta_j$'s and obtain
\[
 J = \pm \prod_{k=1}^{d-1} |\mathbf{\Theta}_k^+| = \pm \prod_{k=1}^{d-1} |c_{k,k}| = \pm |\{v_{k,l}\}|
\]
Next we observe that the matrix $\{v_{k,l}\}$ is actually the product of the square matrices $\{y_j^{(k)}\}$ and $\{(\partial \Theta_j / \partial \omega_l)\cdot \Theta_0\}$. Therefore we have 
\[
 J = \pm \left|\{y_j^{(k)}\}\right|\cdot \left|\left\{\frac{\partial \Theta_j}{\partial \omega_l}\cdot \Theta_0 \right\}\right| = |A_0 A_1 \cdots A_{d-1}| \cdot \left|\left\{\frac{\partial \Theta_0}{\partial \omega_l}\cdot \Theta_j \right\}\right|
\]
Here we use the fact $\Theta_j \cdot \Theta_0 = 0$. This immediately gives  
\[
  d\mathbf{x}^{(0)} \cdots d\mathbf{x}^{(d)} = |A_0 A_1 \cdots A_{d-1}| \cdot \left|\left\{\frac{\partial \Theta_0}{\partial \omega_l}\cdot \Theta_j \right\}\right| d\omega d\mathbf{y}^{(0)} \cdots d\mathbf{y}^{(d)}.
\]
Finally we recall that $\Theta_j$ are a set of orthonormal basis of the tangent space of $\mathbb{S}^{d-1}$ at the point $\Theta_0 \in \mathbb{S}^{d-1}$, and conclude
\[
 d\mathbf{x}^{(0)} \cdots d\mathbf{x}^{(d)} = |A_0 A_1 \cdots A_{d-1}| d\mathbb{S}^{d-1}(\omega) d\mathbf{y}^{(0)} \cdots d\mathbf{y}^{(d)}.
\]
\end{proof}

\subsection{Proof of geometric inequality}

 We divide all the reciprocal simplexes into a few groups. More precisely, we define  
\[
 \Sigma_h = \left\{(X_0, X_1, \cdots, X_d)\in \Sigma \cap \Omega_\ast^{d+1}: h\leq \hbox{dist} (X_d, X_0 X_1 \cdots X_{d-1}) < 2h\right\}.
\]
for $h \in \{r, r/2, r/4, \cdots\}$. Thus we have 
\[
 \int_{\Sigma \cap \Omega_\ast^{d+1}} \frac{dX_0 dX_1 \cdots dX_d}{|X_0 X_1 \cdots X_d|} \leq \sum_{h} \int_{\Sigma_h}  \frac{dX_0 dX_1 \cdots dX_d}{|X_0 X_1 \cdots X_d|}.
\]
We then apply the change of variable $\mathbf{T}$ defined in Subsection \ref{sec: change of variables} and rewrite the integral above in the form of 
\begin{align*}
 \int_{\Sigma_h}  \frac{dX_0 dX_1 \cdots dX_d}{|X_0 X_1 \cdots X_d|} &= \int_{\mathbf{T} \Sigma_h} \frac{|X_0 X_1 \cdots X_{d-1}| d \mathbb{S}^{d-1} (\omega) dy^{(0)} dy^{(1)} \cdots dy^{(d)}}{|X_0 X_1 \cdots X_d|}\\
 & \lesssim_1 \frac{1}{h} \int_{\mathbf{T} \Sigma_h} d \mathbb{S}^{d-1} (\omega) dy^{(0)} dy^{(1)} \cdots dy^{(d)}
\end{align*}
We first fix\footnote{We slightly abuse the notation to let $\omega_0$ represent the the unit vector $\Theta_0(\omega) \in \mathbb{S}^{d-1}$.} $\omega_0 \in \mathbb{S}^{d-1}$ and $\mathbf{y}^{(0)} \in \Omega_\ast$, if $(X_0, X_1, X_2, \cdots X_d) \in \Sigma_h$, then $X_1, X_2, \cdots, X_{d-1}$ are all contained in the $(d-1)$-dimensional hyperplane which contains the point $X_0$ (uniquely determined by $\mathbf{y}^{0}$) and is orthogonal to the vector $\omega$. Lemma \ref{intersection sphere shell} implies that $X_1, X_2, \cdots, X_{d-1}$ are contained in a $(d-1)$-dimensional sphere shell (or a sphere) of outer radius $r_\ast$ and thickness $w_\ast$. Since $\mathbf{y}^{(k)}$'s are the coordinates of $X_k$ in the hyperplane, we obtain that the set 
\[
 \left\{\mathbf{y}^{(k)}\in \Rm^{d-1}: (\omega_0, \mathbf{y}^{(0)}, \mathbf{y}^{(k)}) \in \mathbf{P}_{\omega, 0, k} \mathbf{T} \Sigma_h\right\} 
\]
is contained in a sphere shell (or a sphere) of the same size in $\Rm^{d-1}$. Here $\mathbf{P}_{\omega, 0,k}$ is the natural projection from $\mathbb{S}^{d-1} \times \Rm^d \times (\Rm^{d-1})^{d-1} \times \Rm^d$ onto $\mathbb{S}^{d-1} \times \Rm^d \times \Rm^{d-1}$ defined by 
\[
 \mathbf{P}_{\omega, 0,k} (\omega, \mathbf{y}^{(0)}, \mathbf{y}^{(1)}, \cdots, \mathbf{y}^{(d)}) = (\omega, \mathbf{y}^{(0)}, \mathbf{y}^{(k)})
\]
The projection operators $\mathbf{P}_{\omega, 0, d}$ and $\mathbf{P}_{\omega,0}$ are defined in a similar way and will be used later in the proof. The shape and size of the set immediately gives 
\begin{equation} \label{size of yk}
 \left|\left\{\mathbf{y}^{(k)}\in \Rm^{d-1}: (\omega_0, \mathbf{y}^{(0)}, \mathbf{y}^{(k)}) \in \mathbf{P}_{\omega, 0, k} \mathbf{T} \Sigma_h\right\}\right| \lesssim_1 r_\ast^{d-2} w_\ast \lesssim_1 r^{d-2} w. 
\end{equation}
We then recall our assumption $h\leq \hbox{dist} (X_d, X_0 X_1 \cdots X_{d-1}) < 2h$ and the choice of $\omega$ to obtain that $h \leq y_0^{(d)} < 2h$. A similar argument as above then gives 
\begin{equation} \label{size of yd}
 \left|\left\{\mathbf{y}^{(d)}\in \Rm^{d-1}: (\omega_0, \mathbf{y}^{(0)}, \mathbf{y}^{(d)}) \in \mathbf{P}_{\omega, 0, d} \mathbf{T} \Sigma_h\right\}\right| \lesssim_1 r^{d-2} wh. 
\end{equation} 
We then utilize the upper bounds \eqref{size of yk}, \eqref{size of yd} in the integral above and obtain 
\begin{align}
 \int_{\Sigma_h}  \frac{dX_0 dX_1 \cdots dX_d}{|X_0 X_1 \cdots X_d|}  & \lesssim_1 \frac{1}{h} \int_{\mathbf{P}_{\omega,0} \mathbf{T} \Sigma_h} (r^{d-2} w)^{d-1} (r^{d-2} wh) d \mathbb{S}^{d-1} (\omega) d\mathbf{y}^{(0)} \nonumber\\
 & \lesssim_1  r^{d(d-2)} w^d \int_{\mathbf{P}_{\omega,0} \mathbf{T} \Sigma_h} d \mathbb{S}^{d-1} (\omega) d\mathbf{y}^{(0)}.  \label{integral upper bound 2}
\end{align}
Next we fix $\omega \in \mathbb{S}^{d-1}$ and consider the possible location of $\mathbf{y}^{(0)}$, or equivalently $X_0$ in $\Rm^d$. We apply Lemma \ref{location of reciprocal} and obtain that if $(\omega, \mathbf{y}^{(0)}) \in \mathbf{P}_{\omega,0} \mathbf{T} \Sigma_h$, then there exists $k \in\{0,1,\cdots, d\}$, so that 
\[
 \left|\omega\cdot \overrightarrow{X_0 A_k} \right| \lesssim h.
\]
We then apply Lemma \ref{intersection sphere shell} again and obtain 
\[
 \left|\left\{X \in \Omega_\ast: \min_{k \in \{0,1,\cdots,k\}} \left|\omega\cdot \overrightarrow{X A_k} \right| \lesssim h \right\}\right| \lesssim wr^{d-2} h.
\]
This implies 
\[
 \left|\{\mathbf{y}^{(0)} \in \Rm^d: (\omega, \mathbf{y}^{(0)}) \in \mathbf{P}_{\omega,0} \mathbf{T} \Sigma_h\}\right| \lesssim wr^{d-2} h.
\]
This gives the upper bound of the integral in \eqref{integral upper bound 2}
\begin{align*}
  \int_{\Sigma_h}  \frac{dX_0 dX_2 \cdots dX_d}{|X_0 X_1 \cdots X_d|} \lesssim r^{d(d-2)} w^d \int_{\mathbb{S}^{d-1}} (wr^{d-2} h) d\mathbb{S}^{d-1}(\omega)
  \lesssim r^{d^2-d-2} w^{d+1} h. 
\end{align*}
Finally we take a sum for $h \in \{r, r/2, \cdots\}$ and finish the proof. 
\[
  \int_{\Sigma \cap \Omega_\ast^{d+1}} \frac{dX_0 dX_2 \cdots dX_d}{|X_0 X_1 \cdots X_d|} \leq \sum_{h} r^{d^2-d-2} w^{d+1} h \lesssim r^{d^2-d-1} w^{d+1}.
\]

\subsection{A few remarks} \label{sec: optimal upper bound} 

In this subsection we first show that the upper bound given in Theorem \ref{main 1} is optimal and then give a few more remarks on Theorem \ref{main 1}. Let us consider the regular simplex $B_0 B_1 \cdots B_d$ inscribed in the outer surface of the sphere shell and define ($O$ is the center of sphere shell $\Omega_\ast$, the constant $\eps = \eps(d)$ is sufficiently small)
 \[
  \Omega_k = \{X\in \Omega_\ast: \angle XOB_k < \eps; r-\min\{\eps r, w\} <|XO|< r\}. 
 \]
Then the integral 
\[
 J = \int_{\Omega_0 \times \Omega_1 \times \cdots \times \Omega_d} \left(\int_{\Sigma(X_0, \cdots, X_d)\cap \Omega_\ast^{d+1}} \frac{dY_0 dY_1\cdots dY_d}{|Y_0 Y_1 \cdots Y_d|}\right) dX_0 dX_1 \cdots dX_d
\]
satisfies
\begin{align*} 
 J & \geq \int_{\Omega_0 \times \Omega_1 \times \cdots \times \Omega_d} \left(\int_{\Omega_0 \times \Omega_1 \times \cdots \times \Omega_d} \frac{\psi(X_0, \cdots, X_d, Y_0, \cdots, Y_d) dY_0 dY_1\cdots dY_d}{|Y_0 Y_1 \cdots Y_d|}\right) dX_0 dX_1 \cdots dX_d\\
 & \gtrsim_d  \frac{1}{r^d}\int_{\Omega_0 \times \Omega_1 \times \cdots \times \Omega_d \times \Omega_0 \times \Omega_1 \times \cdots \times \Omega_d} \psi(X_0, \cdots, X_{d}; Y_0, \cdots, Y_d) dX_0 \cdots dX_{d} dY_0 \cdots dY_{d}.
\end{align*} 
Here $\psi(X_0, \cdots, X_{d}; Y_0, \cdots Y_d)$ is the ``characteristic function'' of reciprocal simplexes, i.e.
\[
 \psi(X_0, \cdots, X_{d}; Y_0, \cdots Y_d) = \left\{\begin{array}{ll} 1, & X_0 X_1 \cdots X_{d} \; \hbox{and}\; Y_{0} Y_{1} \cdots Y_{d}\; \hbox{are reciprocal};\\
 0, & \hbox{otherwise.}\end{array}\right.
\]
We may switch the variables $X_0 \leftrightarrow Y_0, X_1 \leftrightarrow Y_1, \cdots$ and define a family of transformations 
\[
 \mathbf{S}_{j_0 j_1 \cdots j_d}: \Omega_0 \times \Omega_1 \times \cdots \times \Omega_d \times \Omega_0 \times \Omega_1 \times \cdots \times \Omega_d \rightarrow \Omega_0 \times \Omega_1 \times \cdots \times \Omega_d \times \Omega_0 \times \Omega_1 \times \cdots \times \Omega_d
\]
Here $j_0, j_1, \cdots, j_d \in \{0,1\}$; the transformation $\mathbf{S}_{j_0 j_1 \cdots j_d}(X_0, X_1, \cdots, X_d; Y_0, Y_1, \cdots, Y_d)$ is defined in the following way: we switch $X_k$ and $Y_k$ if and only if $j_k = 1$. For example, 
\[
 \mathbf{S}_{010\cdots01}(X_0, X_1, \cdots, X_d; Y_0, Y_1, \cdots, Y_d) = (X_0,Y_1,X_2,\cdots, X_{d-1},Y_d; Y_0, X_1, Y_2, \cdots, Y_{d-1}, X_d).
\]
We then switch the variables and obtain 
\begin{align*}
 J \gtrsim_d \frac{1}{r^d} \sum_{j_k \in \{0,1\}} \int_{(\Omega_0 \times \Omega_1 \times \cdots \times \Omega_d)^2} \psi(\mathbf{S}_{j_0 j_1 \cdots j_d}(X_0, \cdots, X_{d}; Y_0, \cdots Y_d)) dX_0 \cdots dX_{d} dY_0 \cdots dY_{d}.
\end{align*}
We then make an observation: Because any simplex whose vertices include $X_k, Y_k$ has a small volume ($\lesssim_d \eps r^d$), which is much smaller than the volume of simplexes $X_0 X_1 \cdots X_d$ or $Y_0 Y_1\cdots Y_d$, it immediately follows that if we split $X_0, \cdots, X_d, Y_0, \cdots, Y_d$ into two reciprocal simplexes, then each simplex must contain $d+1$ points with different subscripts, although these points may be a mixture of $X$'s and $Y$'s. Therefore we have 
\[
 \sum_{j_k \in \{0,1\}}  \psi(\mathbf{S}_{j_0 j_1 \cdots j_d}(X_0, \cdots, X_{d}; Y_0, \cdots Y_d)) \geq 2.
\]
Thus 
\[
 J \gtrsim_d  \frac{1}{r^d} \int_{(\Omega_0 \times \Omega_1 \times \cdots \times \Omega_d)^2} dX_0 \cdots dX_{d} dY_0 \cdots dY_{d} \gtrsim_d \frac{|\Omega_0|^2 |\Omega_1|^2\cdots |\Omega_d|^2}{r^d}.
\]
We recall the definition of $J$ and obtain 
\[
 \sup_{(X_0, \cdots, X_d)\in \Omega_0 \times \cdots \times \Omega_d} \int_{\Sigma(X_0, \cdots, X_d)\cap \Omega_\ast^{d+1}} \frac{dY_0 dY_1\cdots dY_d}{|Y_0 Y_1 \cdots Y_d|} \gtrsim_d \frac{|\Omega_0| |\Omega_1| \cdots |\Omega_d|}{r^d} \gtrsim_d w^{d+1} r^{d^2-d-1}.
\]
This shows that the upper bound $w^{d+1} r^{d^2-d-1}$ is optimal up to a constant depending on $d$. Before we finish this section we give a few more remarks. 

\begin{remark}
 The proof given in this section seems to be easier than the argument we used in \cite{inequality26} for reciprocal triangles. But the argument here does not apply to the case of reciprocal triangles directly because the inequality $r_\ast^{d-2} w_\ast \lesssim r^{d-2}w$ in \eqref{size of yk} does not hold when $d=2$ although we still have $r_\ast w_\ast \lesssim rw$. 
\end{remark}

\begin{remark}
A review of our proof shows that Theorem \ref{main 1} also holds if we use a weaker version of reciprocal simplexes instead. More precisely we call two polytopes $A_1 A_2 \cdots A_{d+1}$ and $A_{d+2} A_{d+3} \cdots A_{2d+2}$ in $\Rm^d$ weakly reciprocal to each other, if and only if the product of their volumes satisfies 
 \[
  |A_1 A_2 \cdots A_{d+1}|\cdot |A_{d+2} A_{d+3} \cdots A_{2d+2}| \geq \gamma \max |A_{j_1} A_{j_2} \cdots A_{j_{d+1}}|\cdot |A_{j_{d+2}} A_{j_{d+3}} \cdots A_{j_{2d+2}}|.
 \]
 Here $\gamma\in (0,1)$ is a constant. The maximum is again taken for all permutations $\{j_1,j_2,\cdots,j_{2d+2}\}$ of the first $2d+2$ positive integers $1,2,\cdots, 2d+2$. 
\end{remark}

\section{Decay Estimates of Adjoint Radon Transforms}

In this section we prove the decay estimates of the adjoint Radon transform $\mathcal{R}^*$ given in Theorem \ref{main 2}, as well as Corollary \ref{global layer like Lp}. The first and most important step is to prove part (a) of Theorem \ref{main 2}, i.e. decay estimates with localized data away from the origin.

\subsection{Localized data away from the origin}

The general idea is the same as the three dimensional case. First of all, it suffices to consider continuous functions $G$ supported in $[a,b] \times \mathbb{S}^{d-1}$, by a standard approximation technique. In this case $G$ is actually uniformly continuous. We introduce 
\[
 v(x,t) = \int_{\mathbb{S}^{d-1}} G(x\cdot \omega + t, \omega) d\omega
\]
By the uniform continuity of $G$, we have the uniform converge $v(x,t) \rightrightarrows \mathcal{R}^* G (x)$ as $t\rightarrow 0$. This gives the uniform converge 
\[
 \frac{1}{\delta}\int_0^\delta v(x,t) dt \rightrightarrows \mathcal{R}^* G (x), \qquad \hbox{as}\; \delta \rightarrow 0^+.
\]
As a result, Fatou's lemma gives
\[
 \int_{|x|>R} |\mathcal{R}^* G(x)|^{2d} dx \leq \liminf_{\delta\rightarrow 0^+} \int_{|x|>R} \left|\frac{1}{\delta}\int_0^\delta v(x,t) dt\right|^{2d} dx. 
\]
For convenience we define $J_\delta$ to be the integral in the right hand side of inequality above: 
\[
 J_\delta = \int_{|x|>R} \left|\frac{1}{\delta}\int_0^\delta v(x,t) dt\right|^{2d} dx. 
\]
We may apply the change of variable $s=x\cdot \omega + t$ and write 
\[
\frac{1}{\delta}\int_0^\delta v(x,t) dt = \frac{1}{\delta}\int_{\mathbb{S}^{d-1}} \int_a^b G(s,\omega) \chi(s,\omega,x) ds d\omega 
\]
Here $\chi(s,\omega,x)$ is the characteristic function of the region $\{(s,\omega,x): s-\delta < x\cdot \omega < s\}$. We insert this into $J_\delta$ and obtain
\[
 J_\delta  \leq \frac{1}{\delta^{2d}} \int_{\Rm^d} \int_{(I\times \mathbb{S}^{d-1})^{2d}} \left(\prod_{k=1}^{2d} |G(s_k, \omega_k)| \chi_{A_k} (x) \right) (d\omega ds)^{2d} dx.
\]
Here $I = [a,b]$ and $(d\omega ds)^{2d} = d\omega_1 ds_1 d\omega_2 ds_2\cdots \omega_{2d} ds_{2d}$. In addition, $\chi_{A_k}(x)$ is the characteristic function of 
\[
 A_k = \{x\in \Rm^d: |x|>R, s-\delta < x\cdot \omega < s\}.
\]
We then introduce the reciprocal groups of unit vectors. We call two groups of unit vectors $(\omega_1, \omega_2, \cdots, \omega_d), (\omega_{d+1}, \omega_{d+2}, \cdots, \omega_{2d}) \in (\mathbb{S}^{d-1})^d$ reciprocal to each other if 
\[
 \left| \omega_1, \omega_2, \cdots, \omega_d \right| \cdot \left| \omega_{d+1}, \omega_{d+2}, \cdots, \omega_{2d} \right| = \max_{j_1, j_2, \cdots, j_{2d}} \left| \omega_{j_1}, \omega_{j_2}, \cdots, \omega_{j_d} \right| \cdot \left| \omega_{j_{d+1}}, \omega_{j_{d+2}}, \cdots, \omega_{j_{2d}} \right|.
\]
The maximum is taken for all permutations $j_1, j_2, \cdots, j_{2d}$ of $1,2,\cdots, 2d$. Here the notation $|\omega_1, \omega_2, \cdots, \omega_d|$ denotes the volume of the parallel polyhedron spanned by the vectors $\omega_1, \omega_2, \cdots, \omega_d$. It is clear that given $\omega_1, \omega_2, \cdots, \omega_{2d} \in \mathbb{S}^{d-1}$, we can always divide them into two reciprocal groups of $d$ vectors. Therefore by rotating the variables we only need to consider the integral in the region $\Sigma\subset (\mathbb{S}^{d-1})^{2d}$ where $(\omega_1, \omega_2, \cdots, \omega_d)$ and $(\omega_{d+1}, \omega_{d+2}, \cdots, \omega_{2d})$ are reciprocal. We then apply Cauchy-Schwarz inequality ($A_{j_1, j_2, \cdots, j_{m}} = A_{j_1}\cap A_{j_2} \cap \cdots \cap A_{j_m}$)
\begin{align*}
 J_\delta & \lesssim_d \frac{1}{\delta^{2d}} \int_{|x|>R} \int_{\Sigma \times I^{2d}} \left(\prod_{k=1}^{2d} |G(s_k, \omega_k)| \chi_{A_k} (x) \right) (d\omega ds)^{2d} dx\\
 & \lesssim_d \frac{1}{\delta^{2d}} \int_{\Rm^d} \int_{\Sigma \times I^{2d}} \frac{|A_{d+1, d+2, \cdots, 2d}|}{|A_{1,2,\cdots ,d}|} |G(s_1,\omega_1)|^2  \cdots  |G(s_d,\omega_d)|^2 \chi_{A_{1,2,\cdots,2d}}(x) (d\omega ds)^{2d} dx\\
 & \quad + \frac{1}{\delta^{2d}} \int_{\Rm^d} \int_{\Sigma \times I^{2d}} \frac{|A_{1,2,\cdots ,d}|}{|A_{d+1, d+2, \cdots, 2d}|} |G(s_{d+1},\omega_{d+1})|^2 \cdots |G(s_{2d},\omega_{2d})|^2\chi_{A_{1,2,\cdots,2d}}(x) (d\omega ds)^{2d} dx\\
 & \lesssim_d \frac{1}{\delta^{2d}} \int_{\Rm^d} \int_{\Sigma \times I^{2d}} \frac{|A_{d+1, d+2, \cdots, 2d}|}{|A_{1,2,\cdots ,d}|} |G(s_1,\omega_1)|^2  \cdots  |G(s_d,\omega_d)|^2 \chi_{A_{1,2,\cdots,2d}}(x) (d\omega ds)^{2d} dx.
\end{align*}
Thus 
\begin{align}
 J_\delta & \lesssim_d \int_{(\mathbb{S}^{d-1})^d\times I^d} J(s_{1,2,\dots,d}, \omega_{1,2,\cdots, d}) |G(s_1,\omega_1)|^2  \cdots  |G(s_d,\omega_d)|^2  ds_{1,2,\cdots, d} d\omega_{1,2,\cdots, d} \nonumber\\
 & \lesssim_d \|G\|_{L^2}^{2d} \sup_{s_{1,2,\cdots, d}\in I^d, \omega_{1,2,\cdots, d}\in (\mathbb{S}^{d-1})^d} J(s_{1,2,\dots,d}, \omega_{1,2,\cdots, d}). \label{JdeltaJ}
\end{align} 
Here we use the notations $s_{j_1, j_2, \cdots, j_m} = (s_{j_1}, s_{j_2}, \cdots, s_{j_m})$, $d s_{j_1, j_2, \cdots, j_m} = ds_{j_1} ds_{j_2} \cdots ds_{j_m}$, and 
\begin{align*}
 J(s_{1,2,\dots,d}, \omega_{1,2,\cdots, d}) & = \frac{1}{\delta^{2d}} \int_{\Rm^d} \int_{\Sigma(\omega_{1,2,\cdots,d})\times I^d} \frac{|A_{d+1, d+2, \cdots, 2d}|}{|A_{1,2,\cdots ,d}|}\chi_{A_{1,2,\cdots,2d}}(x) (ds d\omega)^d dx;\\
 \Sigma(\omega_{1,2,\cdots,d}) & = \{(\omega_{d+1}, \cdots, \omega_{2d}): (\omega_1, \cdots, \omega_d)\, \hbox{and}\, (\omega_{d+1}, \cdots, \omega_{2d})\, \hbox{are reciprocal}\}
\end{align*}
For convenience we slightly abuse the notation above $(ds d\omega)^d = ds_{d+1} d\omega_{d+1} \cdots ds_{2d} d\omega_{2d}$. The remaining work is to find an upper bound of $J(s_{1,2,\dots,d}, \omega_{1,2,\cdots, d}) $. We have 
\begin{align*}
 J &\leq \frac{1}{\delta^{2d}|A_{1,2,\cdots,d}|} \int_{A_{1,2,\cdots, d}} \int_{\Sigma(\omega_{1,2,\cdots,d})\times I^d}        |A_{d+1, d+2, \cdots, 2d}| \chi_{A_{d+1, \cdots, 2d}} (x) (ds d\omega)^d dx \\
 & \leq \frac{1}{\delta^{2d}} \sup_{x\in A_{1,2,\cdots,d}} \int_{\Sigma(\omega_{1,2,\cdots,d})\times I^d}    |A_{d+1, d+2, \cdots, 2d}| \chi_{A_{d+1, \cdots, 2d}} (x) (ds d\omega)^d.
\end{align*}
Given $x \in B_R^c \doteq \{y: |y|>R\}$, we define 
\[
 \Omega_\delta (x) = \{\omega\in \mathbb{S}^{d-1}: \exists s\in I, s-\delta < x\cdot \omega<s\} = \{\omega\in \mathbb{S}^{d-1}: a-\delta < x\cdot \omega<b\},
\]
and further discuss the upper bound of $\sup J$
\begin{align*}
 \sup J & \leq \sup_{s_{1,\cdots, d}, \omega_{1,\cdots,d}} \left(\sup_{x\in A_{1,2,\cdots,d}} \delta^{-2d}\int_{\Sigma(\omega_{1,2,\cdots,d})\times I^d}    |A_{d+1, d+2, \cdots, 2d}| \chi_{A_{d+1, \cdots, 2d}} (x) (ds d\omega)^d\right)\\
 & \leq \sup_{x\in B_R^c; \omega_1, \omega_2, \cdots, \omega_d \in \Omega_\delta(x)} \delta^{-2d} \int_{\left(\Sigma(\omega_{1,2,\cdots,d})\cap (\Omega_\delta(x))^{d}\right)\times I^d}    |A_{d+1, d+2, \cdots, 2d}| \chi_{A_{d+1, \cdots, 2d}} (x) (ds d\omega)^d.
\end{align*}
Geometrically speaking, each set $A_k$ is a ``thin slice'' which is orthogonal to $\omega$ and of thickness $\delta$, possibly with a hole(we delete the ball $|x|<R$). Therefore
\[
 |A_{d+1, d+2, \cdots, 2d}| \leq \frac{\delta^d}{|\omega_{d+1}, \omega_{d+2}, \cdots, \omega_{2d}|}.
\]
In addition, we observe that if $x\in \Rm^d$, $\omega_k \in \mathbb{S}^{d-1}$ are both fixed, then 
\[
 x \in A_k \quad \Rightarrow \quad s_k \in (x\cdot \omega, x\cdot \omega+\delta).
\]
Combining these two facts, we have 
\[
 \sup J \leq \sup_{x\in B_R^c; \omega_1, \omega_2, \cdots, \omega_d \in \Omega_\delta(x)} \int_{\Sigma(\omega_{1,2,\cdots,d})\cap (\Omega_\delta(x))^{d}}  \frac{d\omega_{d+1, \cdots, 2d}}{|\omega_{d+1}, \omega_{d+2}, \cdots, \omega_{2d}|}.
\]
By radial symmetry, it suffices to consider $x = (0,0,\cdots, 0,h)$ with $h > R$. Thus we have 
\[
 \sup J \leq \sup_{h> R} \left(\sup_{\omega_1, \omega_2, \cdots, \omega_d \in \Omega_\delta(h)} \int_{\Sigma(\omega_{1,2,\cdots,d})\cap (\Omega_\delta(h))^{d}}  \frac{d\omega_{d+1, \cdots, 2d}}{|\omega_{d+1}, \omega_{d+2}, \cdots, \omega_{2d}|}\right).
\]
Here 
\begin{align*}
 \Omega_{\delta} (h) & = \left\{\omega = (x_1,x_2,\cdots,x_d) \in \mathbb{S}^{d-1}: \frac{a-\delta}{h}< x_d< \frac{b}{h}\right\}.&  & h\geq b;\\
 \Omega_{\delta} (h) & = \left\{\omega = (x_1,x_2,\cdots,x_d) \in \mathbb{S}^{d-1}: \frac{a-\delta}{h}< x_d \leq 1\right\}, & & h\in (a-\delta,b);\\
 \Omega_{\delta} (h) & = \varnothing. & & h \leq a-\delta.
\end{align*}
The case $h\leq a-\delta$ is trivial. Thus we assume $h>a-\delta$. Next we apply the central projection $\mathbf{P}$ from the upper half sphere $\mathbb{S}_+^{d-1}$ to the hyperplane $x_d = 1$ defined by (the center is the origin $O$ of $\Rm^d$)
\[
 \mathbf{P} (x_1, x_2, \cdots, x_d) = (x_1/x_d, x_2/x_d, \cdots, x_{d-1}/x_d, 1);
\]
and let $Y_k = \mathbf{P} \omega_k$. If we use the notation $\omega_{k,d}$ for the last component of $\omega_{k}$, then we have $\overrightarrow{OY_k} = \omega_{k,d}^{-1} \omega_k$. Therefore the area of $Y_{j_1} Y_{j_2} \cdots Y_{j_d}$ satisfies 
\[
 |Y_{j_1} Y_{j_2} \cdots Y_{j_d}| = d |O Y_{j_1} Y_{j_2} \cdots Y_{j_d}| = \frac{1}{(d-1)!} \cdot \frac{ \left|\omega_{j_1}, \omega_{j_2}, \cdots, \omega_{j_d}\right|}{\omega_{j_1,d} \omega_{j_2,d}\cdots \omega_{j_d,d}}
\]
As a result, for any permutation $j_1, j_2, \cdots, j_{2d}$ of $1,2,\cdots, 2d$ we have 
\begin{align*}
 |Y_1 Y_2 \cdots Y_d|\cdot |Y_{d+1} Y_{d+2}\cdots Y_{2d}| & = \frac{\left|\omega_{1}, \omega_{2}, \cdots, \omega_{d}\right|\cdot \left|\omega_{d+1}, \omega_{d+2}, \cdots, \omega_{2d}\right|}{((d-1)!)^2 \omega_{1,d} \omega_{2,d} \cdots \omega_{2d,d}}  \\
 & \geq \frac{\left|\omega_{j_1}, \omega_{j_2}, \cdots, \omega_{j_d}\right|\cdot \left|\omega_{J_{d+1}}, \omega_{j_{d+2}}, \cdots, \omega_{j_{2d}}\right|}{((d-1)!)^2 \omega_{1,d} \omega_{2,d} \cdots \omega_{2d,d}}\\
 & = |Y_{j_1} Y_{j_2} \cdots Y_{j_d}|\cdot |Y_{j_{d+1}} Y_{j_{d+2}}\cdots Y_{j_{2d}}|
\end{align*}
In other words, $Y_1 Y_2 \cdots Y_d$ and $Y_{d+1} Y_{d+2}\cdots Y_{2d}$ are reciprocal simplexes. We then apply change of variables $(\omega_{d+1}, \omega_{d+2}, \cdots, \omega_{2d}) \rightarrow (Y_{d+1}, Y_{d+2}, \cdots Y_{2d})$ and obtain ($dY_k = \omega_{k,d}^{-d} d\omega_k$)
\begin{align*}
 \sup J &\lesssim_d \sup_{h> \max\{R,a-\delta\}} \left(\sup_{Y_1, \cdots, Y_d \in \Omega_\delta^\ast (h)} \int_{\Sigma(Y_1, \cdots, Y_d)\cap (\Omega_\delta^\ast (h))^{d}}  \frac{(\omega_{d+1,d} \cdots \omega_{2d,d})^{d-1} dY_{d+1} \cdots dY_{2d}}{|Y_{d+1} Y_{d+2} \cdots Y_{2d}|}\right)\\
 & \lesssim_d \sup_{h> \max\{R,a-\delta\}} \frac{b^{d(d-1)}}{h^{d(d-1)}}\left(\sup_{Y_1, \cdots, Y_d \in \Omega_\delta^\ast (h)} \int_{\Sigma(Y_1, \cdots, Y_d)\cap (\Omega_\delta^\ast (h))^{d}}  \frac{dY_{d+1} \cdots dY_{2d}}{|Y_{d+1} Y_{d+2} \cdots Y_{2d}|}\right).
\end{align*}
Here we use the fact $\omega_{k,d} \leq b/h$ for all $w_k \in \Omega_\delta(h)$. The set $\Sigma(Y_1,Y_2, \cdots, Y_d)$ consists of all reciprocal simplexes of $Y_1Y_2 \cdots Y_d$. The region $\Omega_\delta^\ast (h) = \mathbf{P} \Omega_\delta(h)$ is given by
\begin{align*}
 \Omega_\delta^\ast (h) & = \left\{Y\in \Rm^{d-1}: \frac{\sqrt{h^2-b^2}}{b} < |Y| < \frac{\sqrt{h^2-(a-\delta)^2}}{a-\delta}\right\},& & h\geq b; &\\ 
 \Omega_\delta^\ast (h) & = \left\{Y\in \Rm^{d-1}:  |Y| < \frac{\sqrt{h^2-(a-\delta)^2}}{a-\delta}\right\},& &h\in (a-\delta,b).&
\end{align*}
Now we may apply Theorem \ref{main 2} (in the $(d-1)$-dimensional case) to give a upper bound of the integral above and obtain
\[
 \sup J \lesssim_d \sup_{h> \max\{R,a-\delta\}} g(a-\delta, b, h)
\]
with 
\[
 g(a,b,h) = \left\{\begin{array}{ll} \frac{(h-a)^{(d^2-4d+1)/2}(b-a)^d}{b h^{(d^2-2d-1)/2}} & h\geq b; \\ \left(\frac{h-a}{b}\right)^{(d-1)^2/2}, & h \in (a,b). \end{array}\right.
\]
We observe that if $a,b$ are fixed so that $1<b/a\leq \frac{d^2-2d-1}{2(d-1)}$, then $g(a,b,h)$ is an increasing function of $h$ for $h \in \left(a, \frac{(d^2-2d-1)a}{2(d-1)}\right]$ and a decreasing function of $h$ for $h \in \left[\frac{(d^2-2d-1)a}{2(d-1)},+\infty\right)$. Therefore if $\delta \ll a$ is sufficiently small, then we have 
\[
 \sup J \lesssim_d \left\{\begin{array}{ll} g(a-\delta,b,R) \lesssim_d \frac{(b-a+\delta)^d}{(a-\delta) R^{d-1}}, & R\geq \frac{(d^2-2d-1)(a-\delta)}{2(d-1)}; \\ g\left(a-\delta, b, \frac{(d^2-2d-1)(a-\delta)}{2(d-1)}\right) \lesssim_d \left(\frac{b-a+\delta}{a-\delta}\right)^d, & R \in  \left(0, \frac{(d^2-2d-1)(a-\delta)}{2(d-1)}\right). \end{array}\right.
\]
We plugging this upper bound in \eqref{JdeltaJ} and obtain 
\begin{align*}
 J_\delta \lesssim_d \frac{(b-a+\delta)^d}{(a-\delta) (\max\{R, a-\delta\})^{d-1}} \|G\|_{L^2}^{2d}.
\end{align*}
Finally we let $\delta \rightarrow 0^+$ and obtain 
\[
 \int_{|x|>R} \left|\mathcal{R}^* G (x)\right|^{2d} dx \lesssim_d \liminf_{\delta\rightarrow 0^+} J_\delta \lesssim_d \frac{(b-a)^d}{a (\max\{R, a\})^{d-1}} \|G\|_{L^2}^{2d}.
\]

\subsection{Localized Data}

Assume that $G \in L^2(\Rm \times \mathbb{S}^{d-1})$ is supported in $[-b,b]\times \mathbb{S}^{d-1}$. We may fix a constant $\gamma \in \left(\frac{2(d-1)}{d^2-d-1},1\right)$ and split $G$ into a sum 
\[
 G = \sum_{k=0}^\infty G_k. 
\]
The functions $G_k$ here are defined by 
\[
 G_k (s,\omega) = \left\{\begin{array}{ll} G(s, \omega), & s \in [-\gamma^{k} b, -\gamma^{k+1} b) \cup (\gamma^{k+1} b, \gamma^k b];\\
 0, & \hbox{otherwise}.  \end{array}\right.
\]
By symmetry and linearity a similar conclusion to part (a) holds if the function $G$ is supported in $([-b,-a] \cup [a,b])\times \mathbb{S}^{d-1}$ with $1<b/a<\frac{d^2-2d-1}{2(d-1)}$. As a result we have 
\[
  \int_{|x|>R} \left|(\mathcal{R}^* G_k)(x)\right|^{2d} dx \lesssim_d \frac{(\gamma^k b)^{d-1}}{R^{d-1}} \|G_k\|_{L^2}^{2d}.
\]
Therefore 
\[
 \|\mathcal{R}^* G\|_{L^{2d}(\{x: |x|>R\})} \leq \sum_{k=0}^\infty \|\mathcal{R}^* G_k\|_{L^{2d}(\{x: |x|>R\})} \lesssim_d \sum_{k=0}^\infty (\gamma^k b/R)^{\frac{d-1}{2d}} \|G_k\|_{L^2} \lesssim_d (b/R)^{\frac{d-1}{2d}} \|G\|_{L^2}. 
\]

\subsection{Layer-wise $L^p$ estimate}

In this subsection we prove Corollary \ref{global layer like Lp}. It suffices to consider the case $\gamma \in \left(1, \frac{d^2-2d-1}{2(d-1)}\right)$. In fact, if $\gamma = \beta^2$, then we have 
\begin{align*}
 \|\mathcal{R}^* G\|_{L^{2d}(\{x: \gamma^k < |x|\leq \gamma^{k+1}\})}^2 & = \left(\|\mathcal{R}^* G\|_{L^{2d}(\{x: \beta^{2k} < |x| \leq \beta^{2k+1}\})}^{2d} + 
 \|\mathcal{R}^* G\|_{L^{2d}(\{x: \beta^{2k+1} < |x| \leq \beta^{2k+2}\})}^{2d} \right)^{1/d}\\
 & \leq \|\mathcal{R}^* G\|_{L^{2d}(\{x: \beta^{2k} < |x| \leq \beta^{2k+1}\})}^2 + \|\mathbf{R}^* G\|_{L^{2d}(\{x: \beta^{2k+1} < |x| \leq \beta^{2k+2}\})}^2.
\end{align*}
Therefore if Corollary \ref{global layer like Lp} holds for $\gamma = \beta$, then it holds for $\gamma = \beta^2$ as well. Now we assume $\gamma \in \left(1, \frac{d^2-2d-1}{2(d-1)}\right)$. Again we rewrite $G$ as a sum
\begin{align*}
 &G = \sum_{k=-\infty}^\infty G_k;& &G_k (s,\omega) = \left\{\begin{array}{ll} G(s, \omega), & s \in [-\gamma^{k+1}, -\gamma^{k}) \cup (\gamma^{k}, \gamma^{k+1}];\\
 0, & \hbox{otherwise}.  \end{array}\right.&
\end{align*}
It is clear that 
\[
 \sum_{k=-\infty}^\infty \|G_k\|_{L^2}^2 = \|G\|_{L^2}^2. 
\]
We observe that 
\begin{align*}
 &(\mathcal{R}^* G_k)(x) = 0, \; |x|\leq \gamma^k;& &\|\mathcal{R}^* G_k\|_{L^{2d}(\{x: |x|>\gamma^j\})} \lesssim_d \gamma^{\frac{(d-1)(k-j)}{2d}} \|G_k\|_{L^2},\; j\geq k.&
\end{align*}
Therefore we have 
\[ 
 \|\mathcal{R}^* G\|_{L^{2d}(\{x: \gamma^j < |x| \leq \gamma^{j+1}\})} \leq \sum_{k=-\infty}^j \|\mathcal{R}^* G_k\|_{L^{2d}(\{x: |x|>\gamma^j\})} \lesssim_d \sum_{k=-\infty}^j \gamma^{\frac{(d-1)(k-j)}{2d}} \|G_k\|_{L^2}.
\]
Next we fix a constant $\gamma_1 \in (1,\gamma)$ and obtain 
\[
 \|\mathcal{R}^* G\|_{L^{2d}(\{x: \gamma^j < |x| \leq \gamma^{j+1}\})}^2 \lesssim_{d,\gamma} \sum_{k=-\infty}^j \gamma_1^{\frac{(d-1)(k-j)}{d}} \|G_k\|_{L^2}^2. 
\]
Finally we take a sum in $j$ and obtain 
\begin{align*}
 \sum_{j=-\infty}^\infty \|\mathcal{R}^* G\|_{L^{2d}(\{x: \gamma^j < |x| \leq \gamma^{j+1}\})}^2 & \lesssim_{d,\gamma} \sum_{j=-\infty}^\infty \sum_{k=-\infty}^j \gamma_1^{\frac{(d-1)(k-j)}{d}} \|G_k\|_{L^2}^2\\
 & \lesssim_{d,\gamma} \sum_{k=-\infty}^\infty \sum_{j=k}^\infty \gamma_1^{\frac{(d-1)(k-j)}{d}} \|G_k\|_{L^2}^2\\
 & \lesssim_{d,\gamma} \sum_{k=-\infty}^\infty \|G_k\|_{L^2}^2 = \|G\|_{L^2}^2.
\end{align*} 

\section{Application on Non-radiative solutions}

As an application of the decay estimates of the adjoint Radon transform, we give decay estimates of non-radiative solutions to the free wave equation in $\Rm^5$. Similar decay estimates can be given in higher odd dimensional space but the argument might be a little more complicated. In dimension $5$, the free wave $u$ can be given explicitly by its radiation profile $G \in L^2(\Rm \times \mathbb{S}^{4})$
\[
 u(x,t) = \mathbf{T}_w G = \int_{\mathbb{S}^4} G'(x\cdot \omega+t, \omega) d\omega.
\]
Here $G'$ is the partial derivative with respect to the first variable $s$. An application of Theorem \ref{main 2} immediately gives decay estimate if $G \in \dot{H}^1 (\Rm \times \mathbb{S}^4)$ is compactly supported. Here the $\dot{H}^1$ space is the completion of the space of $C_0^\infty$ functions equipped with the norm 
\[
 \|G\|_{\dot{H}^1(\Rm \times \mathbb{S}^4)} = \|G'(s,\omega)\|_{L^2(\Rm \times \mathbb{S}^4)}. 
\] 
In order to obtain decay estimate for compactly supported functions in $L^2(\Rm \times \mathbb{S}^4)$, we conduct an interpolation argument. 

\paragraph{Modified cut-off operator} We first consider the modified cut-off operator 
\[
 (\mathbf{P} G)(s, \omega) = \varphi(s) \left[G(s,\omega) - \frac{1}{\|\varphi\|_{L^1}}\int_{\Rm} \varphi(s') G(s',\omega) ds' \right]
\]
Here $\varphi: \Rm \rightarrow [0,1]$ is a smooth even cut-off function satisfying $\varphi(x) = 1$ for $x\in [-1,1]$ and $\varphi(x) = 0$ for $|x|>2$. We claim that $\mathbf{P}$ is a bounded operator from $\dot{H}^1(\Rm \times \mathbb{S}^4)$ to itself. Without loss of generality we assume $G(s,\omega) \in C_0^\infty (\Rm \times \mathbb{S}^4)$. We first fix an $\omega \in \mathbb{S}^4$ and observe 
\[
 \int_{\Rm} \varphi(s) \bar{G}(s,\omega) ds = 0. 
\]
Here for convenience we use the notation 
\[
 \bar{G}(s,\omega) = G(s,\omega) - \frac{1}{\|\varphi\|_{L^1}}\int_{\Rm} \varphi(s') G(s',\omega) ds'.
\]
This immediately gives
\[
 \sup_{s\in [-2,2]} |\bar{G}(s,\omega)| \leq 2 \|\bar{G}'(s,\omega)\|_{L^2([-2,2])}.
\]
Therefore we have 
\begin{align*}
 \int_{\Rm} |\partial_s (\mathbf{P} G)(s, \omega)|^2 ds &\leq 2 \int_{-2}^{2} \left(|\varphi(s) \bar{G}'(s,\omega)|^2 + |\varphi'(s) \bar{G}(s,\omega)|^2 \right)ds \\
 & \lesssim_1 \int_{-2}^{2} |\bar{G}'(s,\omega)|^2 ds \\
 & \lesssim_1 \int_\Rm |G'(s,\omega)|^2 ds
\end{align*}
This completes the proof of boundedness of $\mathbf{P}$. Next we show that $\mathbf{P}$ is also a bounded operator from $\dot{H}^\beta(\Rm\times \mathbb{S}^4)$ for any fixed $\beta \in (-1/2,0)$. Here $\dot{H}^\beta(\Rm\times \mathbb{S}^4)$ is the completion of $C_0^\infty (\Rm \times \mathbb{S}^4)$ equipped with the norm 
\[
 \|u\|_{\dot{H}^\beta (\Rm \times \mathbb{S}^4)} = \left(\iint_{\Rm \times \mathbb{S}^4} |\mathbf{D}^\beta u(s,\omega)|^2 ds d\omega\right)^{1/2}. 
\]
First of all, the multiplication operator $f(s) \rightarrow \varphi(s) f(s)$ is a bounded operator from $\dot{H}^\beta (\Rm)$ to itself thus $G(s,\omega) \rightarrow \varphi(s) G(s,\omega)$ is a bounded operator from $\dot{H}^\beta (\Rm\times \mathbb{S}^4)$. In addition, we have $\varphi(s) \in \dot{H}^\beta(\Rm)$ and 
\[
 \left|\int_{\Rm} \varphi(s') G(s',\omega) ds'\right| \leq \|\varphi\|_{\dot{H}^{-\beta}(\Rm)} \|G(\cdot,\omega)\|_{\dot{H}^\beta(\Rm)}.
\]
Therefore 
\begin{align*}
 \left\|\|\varphi\|_{L^1}^{-1} \varphi(s) \int_{\Rm} \varphi(s') G(s',\omega) ds' \right\|_{\dot{H}^\beta(\Rm \times \mathbb{S}^4)}^2 & = \|\varphi\|_{L^1}^{-2} \int_{\mathbb{S}^4} \left|\int_{\Rm} \varphi(s') G(s',\omega) ds'\right|^2 \|\varphi\|_{\dot{H}^\beta(\Rm)}^2 d\omega\\
 & \lesssim \int_{\mathbb{S}^4} \|G(\cdot,\omega)\|_{\dot{H}^\beta(\Rm)}^2 d\omega.
\end{align*}
This gives the boundedness of $\mathbf{P}$ in the space $\dot{H}^\beta (\Rm\times \mathbb{S}^4)$. 

\paragraph{Interpolation argument} We then give a decay estimate of $\mathbf{T}_w \mathbf{P} G$ via an interpolation argument. We first consider two endpoints, i.e. the cases $G\in \dot{H}^1(\Rm \times \mathbb{S}^4)$ and $G \in \dot{H}^{-1/6}(\Rm \times \mathbb{S}^4)$. If $G\in \dot{H}^1(\Rm \times \mathbb{S}^4)$, then $\mathbf{P} G \in \dot{H}^1(\Rm \times \mathbb{S}^4)$ is supported in $[-2,2]\times \mathbb{S}^4$. As a result, given $t\in \Rm$, the time translation version $(\mathbf{P}G)_{(t)}(s,\omega)\doteq (\mathbf{P} G) (s+t, \omega)$ is supported in $[-2-t, 2-t]\times \mathbb{S}^4$. Next we observe the basic identity 
\[ 
 (\mathbf{T}_w \mathbf{P} G)(\cdot,t) = \mathcal{R}^* (\mathbf{P} G)'_{(t)}, \qquad \forall t \in \Rm.
\]
We apply Corollary \ref{translation narrow} and obtain that for any $t\in \Rm$:
\begin{align*}
 \|(\mathbf{T}_w \mathbf{P} G)(\cdot,t)\|_{L^{10} (\{x: |x|>R\})} & = \|\mathcal{R}^* (\mathbf{P} G)'_{(t)}\|_{L^{10}(\{x: |x|>R\})}\\
 & \lesssim_1 R^{-2/5} \|\mathbf{P} G\|_{\dot{H}^1(\Rm \times \mathbb{S}^4)} \lesssim_1 R^{-2/5} \|G\|_{\dot{H}^1(\Rm \times \mathbb{S}^4)}.
\end{align*}
Similarly part (a) of Theorem \ref{main 2} gives the convergence
\[
 \|(\mathbf{T}_w \mathbf{P} G)(\cdot,t)\|_{L^{10} (\Rm^5)} \lesssim \left(\frac{|t|+2}{|t|-2}-1\right)^{5} \|G\|_{\dot{H}^1(\Rm \times \mathbb{S}^4)} \rightarrow 0,
\]
as $t\rightarrow \pm \infty$. In addition, we observe that $(\mathbf{P}G)'_{(t)} \in C(\Rm_t; L^2(\Rm \times \mathbb{S}^4))$ is continuous in $t$. Since $\mathcal{R}^*$ is a bounded operator from $L^2(\Rm \times \mathbb{S}^4)$ to $L^{10}(\Rm^5)$, we obtain that $\mathbf{T}_w \mathbf{P} G \in C(\Rm_t; L^{10} (\Rm^5))$.  In summary, we have 
\begin{equation} \label{endpoint 1}
 \|\mathbf{T}_w \mathbf{P} G\|_{L_0^\infty(\Rm_t; L^{10}(\{x: |x|>R\}))} \lesssim_1 R^{-2/5} \|G\|_{\dot{H}^1(\Rm \times \mathbb{S}^4)}. 
\end{equation}
Here the space $L_0^\infty$ is a subspace of the regular $L^\infty$ space defined in the following way: We first consider the space of all simple functions so that each simple function is zero except in a set with a finite measure then define the completion of this space in the $L^\infty$ space to be $L^\infty_0$. Next we consider the case $G \in \dot{H}^{-1/6} (\Rm \times \mathbb{S}^4)$. Remark \ref{Hs equivalence} implies that the corresponding initial data $(u_0,u_1)$ of $\mathbf{P} G$ satisfies 
\[
  \|(u_0,u_1)\|_{\dot{H}^{5/6}\times \dot{H}^{-1/6}(\Rm^5)}\lesssim_1 \|\mathbf{P} G\|_{\dot{H}^{-1/6}(\Rm \times \mathbb{S}^4)} \lesssim_1 \|G\|_{\dot{H}^{-1/6}(\Rm \times \mathbb{S}^4)}.
\]
We then apply the Strichartz estimates (see \cite{strichartz}, for instance) to conclude that the corresponding free wave $\mathbf{T}_w \mathbf{P} G$ satisfies 
\begin{equation} \label{endpoint 2}
 \|\mathbf{T}_w \mathbf{P} G\|_{L^2 L^{30/7}(\Rm \times \Rm^5)} \lesssim_1 \|(u_0,u_1)\|_{\dot{H}^{5/6}\times \dot{H}^{-1/6}(\Rm^5)} \lesssim_1 \|G\|_{\dot{H}^{-1/6}(\Rm \times \mathbb{S}^4)}.
\end{equation}
We then combine \eqref{endpoint 1} with \eqref{endpoint 2}, apply the complex interpolation method (see Bergh-L\"{o}fstr\"{o}m \cite{interpolationspace}, for example) to conclude 
\begin{equation} \label{modified cutoff version}
 \|\mathbf{T}_w \mathbf{P} G\|_{L^{7/3} L^{14/3}(\Rm \times \{x\in \Rm^5: |x|>R\})} \lesssim_1 R^{-2/35} \|G\|_{L^2(\Rm \times \mathbb{S}^4)}. 
\end{equation}

\paragraph{Decay of non-radiative solutions} Next we call the definition of $\mathbf{P}$ and write 
\begin{equation} \label{decomposition of Tw varphi}
 \mathbf{T}_w (\varphi(s) G) = \mathbf{T}_w \mathbf{P} G + \mathbf{T}_w \left(\frac{\varphi(s)}{\|\varphi\|_{L^1}}\int_{\Rm} \varphi(s') G(s',\omega) ds'\right)
\end{equation}
Given $\gamma \in \{-1/6, 1\}$, we have 
\begin{align*}
 \left\|\frac{\varphi(s)}{\|\varphi\|_{L^1}}\int_{\Rm} \varphi(s') G(s',\omega) ds'\right\|_{\dot{H}^\gamma(\Rm\times\mathbb{S}^4)}^2 &\lesssim_1 \int_{\mathbb{S}^4} \|\varphi\|_{\dot{H}^\gamma(\Rm)}^2 \left|\int_{\Rm} \varphi(s') G(s',\omega) ds'\right|^2 d\omega\\
 & \lesssim_1  \|\varphi\|_{\dot{H}^\gamma(\Rm)}^2 \|\varphi\|_{L^2(\Rm)}^2 \int_{\mathbb{S}^4} \|G(\cdot,\omega)\|_{L^2(\Rm)}^2 d\omega\\
 & \lesssim_1 \|G\|_{L^2(\Rm \times \mathbb{S}^4)}^2. 
\end{align*}
We may apply Corollary \ref{translation narrow} and Strichartz estimates in the same way as given above to obtain 
\begin{align*} 
 \left\|\mathbf{T}_w\left(\frac{\varphi(s)}{\|\varphi\|_{L^1}}\int_{\Rm} \varphi(s') G(s',\omega) ds'\right)\right\|_{L^\infty L^{10} (\Rm \times \{x: |x|>R\})} & \lesssim_1 R^{-2/5} \|G\|_{L^2(\Rm \times \mathbb{S}^4)}; \\
 \left\|\mathbf{T}_w\left(\frac{\varphi(s)}{\|\varphi\|_{L^1}}\int_{\Rm} \varphi(s') G(s',\omega) ds'\right)\right\|_{L^2 L^{30/7} (\Rm \times \Rm^5)} & \lesssim_1 \|G\|_{L^2(\Rm \times \mathbb{S}^4)}.
\end{align*}
A combination of these two inequality immediately gives 
\[
\left\|\mathbf{T}_w\left(\frac{\varphi(s)}{\|\varphi\|_{L^1}}\int_{\Rm} \varphi(s') G(s',\omega) ds'\right)\right\|_{L^{7/3} L^{14/7} (\Rm \times \{x: |x|>R\})} \lesssim_1  R^{-2/35} \|G\|_{L^2(\Rm \times \mathbb{S}^4)}.
\]
We then combine this inequality with \eqref{modified cutoff version} and \eqref{decomposition of Tw varphi} to obtain 
\[
 \left\|\mathbf{T}_w (\varphi(s) G)\right\|_{L^{7/3} L^{14/7} (\Rm \times \{x: |x|>R\})} \lesssim_1  R^{-2/35} \|G\|_{L^2(\Rm \times \mathbb{S}^4)}. 
\]
Finally we observe that if $G \in L^2(\Rm \times \mathbb{S}^4)$ is supported in $[-1,1]\times \mathbb{S}^4$, then $\varphi(s) G = G$. This finishes the proof. 
\section*{Acknowledgement}
The authors are financially supported by National Natural Science Foundation of China Project 12071339.

\end{document}